\documentclass [11pt]{amsart}
\usepackage[english]{babel}
\usepackage{amsfonts}
\usepackage{amssymb}
\usepackage{amsmath}

\newcounter{remark}
\newcounter{theor}
\newtheorem{thm}{Theorem}[section]
\newtheorem{prop}{Proposition}[section]

\newtheorem{lemma}{Lemma}[section]

\newcommand{\de}{\partial}

\newcommand{\R}{\mathbb{R}}

\newcommand{\mn}{\sqrt{-1}}
\newcommand{\ve}{\varepsilon}
\newcommand{\Ric}{\mathrm{Ric}}
\newcommand{\db}{\overline{\partial}}

\newcommand{\diam}{\mathrm{diam}}

\newcommand{\ti}[1]{\tilde{#1}}

\newcommand{\Ss}{\mathcal{S}}

\newcommand{\ov}[1]{\overline{#1}}

\newcommand{\tr}[2]{\textrm{tr}_{#1} #2}
\newcommand{\F}{\mathcal{F}}

\numberwithin{equation}{section}

\begin{document}

\title[Adiabatic limits of Ricci-flat K\"ahler metrics]{Adiabatic limits of Ricci-flat K\"ahler metrics}
\author{Valentino Tosatti}
\begin{abstract} 
We study adiabatic limits of Ricci-flat K\"ahler metrics on a Calabi-Yau manifold which is the total space of a holomorphic fibration when the volume of the fibers goes to zero. By establishing some new {\em a priori} estimates for the relevant complex Monge-Amp\`ere equation, we show that the Ricci-flat metrics collapse (away from the singular fibers) to a metric on the base of the fibration. This metric has Ricci curvature equal to a Weil-Petersson metric that measures the variation of complex structure of the Calabi-Yau fibers. This generalizes results of Gross-Wilson for $K3$ surfaces to higher dimensions.
\end{abstract}


\address{Department of Mathematics \\ Harvard University \\ Cambridge, MA 02138\\}
  \email{tosatti@math.harvard.edu}
 
\curraddr{Department of Mathematics \\ Columbia University \\ 2990 Broadway, New York, NY 10027}

\maketitle
\section{Introduction}\label{sect1}
In this paper, which is a continuation of \cite{deg}, we study the behaviour of Ricci-flat K\"ahler metrics on a compact Calabi-Yau manifold when the K\"ahler class degenerates to the boundary of the K\"ahler cone. Given a compact Calabi-Yau manifold, a fundamental theorem of Yau \cite{Ya} says that there exists a unique Ricci-flat K\"ahler metric in each K\"ahler class. If we now move the K\"ahler class inside the K\"ahler cone, the corresponding Ricci-flat metrics vary smoothly, as long as the class does not approach the boundary of the K\"ahler cone. The question that we want to address is to understand what happens to the Ricci-flat metrics if the class goes to the boundary of the K\"ahler cone. This question has been raised by Yau and reiterated by others, see for example \cite{yau2, wilson, ctm}. In our previous work \cite{deg} we have studied the case when the limit class has positive volume.
In this paper we consider the case when the limit volume is zero, and we focus on the situation when the Calabi-Yau manifold admits a holomorphic fibration to a lower dimensional space such that the limit class is the pullback of a K\"ahler class from the base.

To state our result, let us introduce some notation. Let $(X,\omega_X)$ be a compact K\"ahler $n$-manifold with $c_1(X)=0$ in $H^2(X,\mathbb{R})$. The condition that $c_1(X)=0$ is equivalent to the requirement that the canonical bundle of $X$ be torsion (see \cite{deg}), and we call $X$ a Calabi-Yau manifold. We assume that there is a holomorphic map $f:X\to Z$ where $(Z,\omega_Z)$ is another compact K\"ahler manifold.
We denote by $Y$ the image of $X$ under $f$, and we assume that $Y$ is an irreducible normal subvariety of $Z$ of dimension $m$ with $0<m<n$, and that the map $f:X\to Y$ has connected fibers. Then $\omega_0=f^*\omega_Z$ is a smooth nonnegative $(1,1)$ form on $X$, whose cohomology class lies on the boundary of the K\"ahler cone. We will also denote by $\omega_Y$ the restriction of $\omega_Z$ to the regular part of $Y$. The map $f:X\to Y$ is an ``algebraic fiber space'' in 
the sense of \cite{laza} and we can find a subvariety $S\subset X$ such that
$Y\backslash f(S)$ is smooth and $f:X\backslash S\to Y\backslash f(S)$ is a smooth submersion ($S$ consists of singular fibers, as well as fibers of dimension strictly larger than $n-m$). Then for any $y\in Y\backslash f(S)$ the fiber $X_y=f^{-1}(y)$ is a smooth $(n-m)$-manifold, equipped with the K\"ahler form $\omega_X|_{X_y}$. Notice that since $f$ is a submersion near $X_y$, we have that $c_1(X_y)=c_1(X)|_{X_y},$
and so the fibers $X_y$ with $y\in Y\backslash f(S)$ are themselves Calabi-Yau.
Yau's theorem \cite{Ya} says that in each K\"ahler class of $X$ there is a unique K\"ahler metric with Ricci curvature identically zero. For each $0<t\leq 1$ we call $\ti{\omega}_t$ the Ricci-flat K\"ahler metric cohomologous to $[\omega_0]+t[\omega_X]$, and we wish to study the behaviour of these metrics when $t$ goes to zero. First of all in \cite{deg} we proved the following

\begin{thm}[Theorem 3.1 of \cite{deg}]
The Ricci-flat metrics $\ti{\omega}_t$ on $X$ have uniformly bounded diameter as $t$ goes to zero.
\end{thm}

The volume of any fiber $X_y$ with respect to  $\ti{\omega}_t$ is comparable to $t^{n-m}$, and the Ricci-flat metrics $\ti{\omega}_t$ 
approach an ``adiabatic limit''. We will show that the metrics $\ti{\omega}_t$ collapse to a K\"ahler metric on $Y\backslash f(S)$.

Our main theorem is the following:

\begin{thm}\label{main2} 
There is a smooth K\"ahler metric $\omega$ on $Y\backslash f(S)$ such that  
the Ricci-flat metrics $\ti{\omega}_t$ when $t$ approaches zero converge to $f^*\omega$ weakly as currents and also in the $C^{1,\beta}_{loc}$ topology of potentials on compact sets of $X\backslash S$, for any $0<\beta<1$. The metric $\omega$ satisfies
\begin{equation}\label{genke}
\Ric(\omega)=\omega_{WP},
\end{equation}
on $Y\backslash f(S)$, where $\omega_{WP}$ is a Weil-Petersson metric measuring the change of complex structures of the fibers. Moreover for any $y\in Y\backslash f(S)$ if we restrict to $X_y$, the metrics $\ti{\omega}_t$ converge to zero in the $C^{1}$ topology of metrics, uniformly as $y$ varies in a compact set of $Y\backslash f(S)$.
\end{thm}

This result generalizes work of Gross and Wilson \cite{gw}, who considered the case when $f:X\to Y=\mathbb{P}^1$ is an elliptically fibered $K3$ surface with $24$ singular fibers of type $I_1$ (see section \ref{examp}). They achieved their result by writing down explicit approximations of the Ricci-flat metrics near the adiabatic limit. A similar approach in higher dimensions seems out of reach at present. Instead, our main technical tool are some new general {\em a priori} estimates for complex Monge-Amp\`ere equations on the total space of a holomorphic fibration.

The limit equation \eqref{genke} has first been explicitly proposed by Song-Tian \cite{st}, where $\omega$ is called a \emph{generalized K\"ahler-Einstein metric}.

Let us briefly explain the meaning of $\omega_{WP}$, referring the reader to section \ref{sect4} for more details. We have already remarked that the smooth fibers $X_y$ of $f$ are themselves Calabi-Yau $(n-m)$-manifolds, polarized by $\omega_X|_{X_y}$. If the canonical bundles of the fibers are actually trivial we get a map from $Y\backslash f(S)$ to the moduli space of polarized Calabi-Yau manifolds and by pulling back the Weil-Petersson metric we get a smooth nonnegative form $\omega_{WP}$ on $Y\backslash f(S)$ (a similar construction goes through in the case when the fibers have torsion canonical bundle). Notice that $\omega_{WP}$ is identically zero precisely when the complex structure of the fibers doesn't change. The appearance of the Weil-Petersson metric in the more general setting of adiabatic limits of constant scalar curvature K\"ahler metrics was observed by Song-Tian \cite{st}, and further studied by Fine \cite{fine} and Stoppa \cite{stoppa}. 

There are two possible situations that we have in mind for our setup: in one case $Y$ is smooth, and then we can just take $Z=Y$. In the second case we take $Z=\mathbb{P}^N$ with $\omega_Z$ the Fubini-Study metric, and then $Y$ is an algebraic variety. A natural class of examples where this situation arises is the following: $X$ is a a projective Calabi-Yau manifold and $L$ is a semiample line bundle over $X$ 
with Iitaka dimension $\kappa(X,L)=m<n$. Then a classical construction of Iitaka (see 2.1.27 in \cite{laza}) gives a holomorphic map $f:X\to\mathbb{P}^N$ exactly as in the setup. Note that if the log Abundance Conjecture holds then every line bundle $L$ with cohomology class on the boundary of the K\"ahler cone and with $\kappa(X,L)=m<n$ is automatically semiample (see \cite{deg}). This is known to hold if $n=2$. 

The organization of the paper is the following. In section \ref{sect2} we set up the problem as a family of degenerating complex Monge-Amp\`ere equations and state our estimates that imply the main result. In section \ref{sect3} we prove {\em a priori} $C^2$ estimates for these equations as well as $C^3$ estimates along the fibers. In section \ref{sect4} we use the estimates to prove our main theorem, and in section \ref{examp} we provide a few examples.\\

{\bf Acknowledgments.}\ I would like to thank my advisor Shing-Tung Yau for suggesting this problem and for constant support. I also thank Chen-Yu Chi, Jian Song, G\'abor Sz\'ekelyhidi and Ben Weinkove for very useful discussions. I was partially supported by a Harvard Merit Fellowship. These results are part of my PhD thesis at Harvard University \cite{tesi}.

\section{Complex Monge-Amp\`ere equations}\label{sect2}
In this section we translate our problem to a family of degenerating complex Monge-Amp\`ere equations, and we state our estimates.\\

First of all let us recall our setup from the introduction: $(X,\omega_X)$ is a compact K\"ahler manifold of complex dimension $n$ with $c_1(X)=0$, or in other words a Calabi-Yau manifold. We have a map holomorphic map $f:X\to Z$, where $(Z,\omega_Z)$ is another compact K\"ahler manifold, with image $Y\subset Z$
and so that $f:X\to Y$ has connected fibers. $Y$ is assumed to be an irreducible normal subvariety of $Z$ of dimension $m$ with $0<m<n$, and we let $\omega_Y$ be the restriction of $\omega_Z$ to the regular part of $Y$.
We also set $\omega_0=f^*\omega_Z$, which is a smooth nonnegative $(1,1)$ form on $X$ whose cohomology class lies on the boundary of the K\"ahler cone. There is a proper subvariety $S\subset X$ such that
$Y\backslash f(S)$ is smooth and $f:X\backslash S\to Y\backslash f(S)$ is a smooth submersion.
Yau's theorem \cite{Ya} says that in each K\"ahler class of $X$ there is a unique K\"ahler metric with Ricci curvature identically zero. For each $0<t\leq 1$ we call $\ti{\omega}_t$ the Ricci-flat K\"ahler metric cohomologous to $[\omega_0]+t[\omega_X]$, and we wish to study the behaviour of these metrics when $t$ goes to zero. On $X$ we have $$\omega_0^k\wedge\omega_X^{n-k}=0,$$
for $m+1\leq k\leq n$, and
\begin{equation}\label{jacob}
\omega_0^m\wedge\omega_X^{n-m}=H \omega_X^n,
\end{equation}
where the smooth non-negative function $H$ vanishes precisely on $S$ and is such that $H^{-\gamma}$ is in $L^1$ for some small $\gamma>0$. This is because $H$ is locally comparable to a sum of squares of holomorphic functions (the minors of the Jacobian of $f$).
In particular it follows that 
$$\int_X\omega_0^m\wedge\omega_X^{n-m}>0.$$
For later purposes we need the following construction. Let $\mathcal{I}$ be the ideal sheaf of $f(S)$ inside $Z$. We cover $Z$ by a finite number of open sets $U_k$ so that on each $U_k$ the ideal $\mathcal{I}$ is generated by holomorphic functions $h_{k,j}$, with $1\leq j\leq N_k$. We then fix $\eta_k$ a partition of unity subordinate to the covering $\{U_k\}$ and we let
\begin{equation}\label{sigma}
\sigma=\sum_{k,j} \eta_k |h_{k,j}|^2,
\end{equation}
if $S\neq\emptyset$ and otherwise we just set $\sigma=1$.
Then $\sigma$ is a smooth nonnegative function on $Z$ with zero locus precisely $f(S)$ and there is a constant $C$ so that on $Z$ we have
\begin{equation}\label{hack3}
\sigma\leq C,\quad 0\leq \mn\de\sigma\wedge\db\sigma\leq C\omega_Z, \quad -C\omega_Z\leq \mn\de\db\sigma \leq C\omega_Z.
\end{equation}
Then for any $y\in Y\backslash f(S)$ we have the inequality
\begin{equation}\label{hack}
\sigma(y)^\lambda\leq C\inf_{X_y}H,
\end{equation}
for some constants $C,\lambda$, and we are free to enlarge $\lambda$ if needed. This is because
both of the function $H$ and $f^*\sigma$ on $X$ are locally comparable to a sum of squares of holomorphic functions and they both have zero set equal to $S$. By taking a log resolution of the ideal sheaf of $S$ inside $X$ and we can assume that $S$ is a divisor with simple normal crossings, and then the holomorphic functions have well defined vanishing orders along the irreducible components of $S$, and \eqref{hack} follows.

In this setting we look at the K\"ahler forms $\omega_t=\omega_0+t\omega_X$ for $0< t\leq 1$, which are cohomologous to the Ricci-flat metrics $\ti{\omega}_t$.
We then define a smooth function $E$ by
$$\Ric(\omega_X)=\mn\de\db E,\quad \int_X e^E\omega_X^n=\int_X\omega_1^n,$$
which is possible thanks to the $\de\db$-lemma.
Then the equation $\Ric(\ti{\omega}_t)=0$ is equivalent to
$$\ti{\omega}_t^n=a_t e^E\omega_X^n,$$
where 
$$a_t=\frac{\int_X \omega_t^n}{\int_X\omega_1^n}.$$
Using the $\de\db$-lemma again, we can find smooth functions $\varphi_t$ for $0<t\leq 1$ so
that $\ti{\omega}_t=\omega_t+\mn\de\db\varphi_t$, $\sup_X\varphi_t=0$ and we have
\begin{equation}\label{degma1}
(\omega_t+\mn\de\db\varphi_t)^n=a_t e^E\omega_X^n.
\end{equation}
Notice that as $t$ approaches zero, the constants $a_t$ behave like
\begin{equation}\label{normconst}
\binom{n}{m}\frac{\int_X\omega_0^m\wedge\omega_X^{n-m}}{\int_X\omega_1^n}t^{n-m}+O(t^{n-m+1}).
\end{equation}
We can then write \eqref{degma1} as
\begin{equation}\label{degma2}
(\omega_t+\mn\de\db\varphi_t)^n=c_t t^{n-m} e^E \omega_X^n,
\end{equation}
where the constant $c_t$ is bounded away from zero and infinity as $t$ goes to zero. Equation \eqref{degma2} has been studied for example in \cite{koltian} where a uniform $L^\infty$ bound on $\varphi_t$ was conjectured. When $m=1$ such a bound can be easily proved using the Moser iteration method (see \cite{st}). The bound in the general case was then proved independently by Demailly and Pali \cite{dpali} and by Eyssidieux, Guedj and Zeriahi \cite{egz2}:
\begin{thm}[\cite{dpali, egz2}]\label{linfty}
There is a constant $C$ that depends only on $X, E, \omega_X, \omega_0$ such that for all $0<t\leq 1$ we have
\begin{equation}\label{linftybound}
\|\varphi_t\|_{L^\infty}\leq C.
\end{equation}
\end{thm}
Our goal is to show higher order estimates for $\varphi_t$ which are uniform on compact sets of $X\backslash S$. 
Notice that since
$$0<\tr{\omega_X}{\ti{\omega}_t}=\tr{\omega_X}{\omega_t}+\Delta_{\omega_X}\varphi_t,$$
and since $\tr{\omega_X}{\omega_t}$ is uniformly bounded, we always have a uniform lower bound for $\Delta_{\omega_X}\varphi_t$.

The following are our main results, and together they imply Theorem \ref{main2}. 

\begin{thm}\label{estimates}
There are constants $A,B,C$ that depend only on the fixed data, so that on $X\backslash S$ and for any $0<t\leq 1$ we have
\begin{equation}\label{c2a}
\frac{t}{Ce^{Ae^{B \sigma^{-\lambda}}}}\omega_X\leq \ti{\omega}_t\leq C e^{Ae^{B\sigma^{-\lambda}}}\omega_X,
\end{equation}
where $\sigma$ is defined by \eqref{sigma}. In particular the Laplacian $\Delta_{\omega_X}\varphi_t$ is bounded uniformly on compact sets of $X\backslash S$, independent of $t$.
\end{thm}

\begin{thm}\label{estimates2}
Given any $y\in Y\backslash f(S)$ denote by $X_y$ the fiber $f^{-1}(y)$, by $\omega_y$ the K\"ahler form
$\omega_X|_{X_y}$ and by $\ti{\omega}_y$ the restriction of the Ricci-flat metric $\ti{\omega}_t|_{X_y}$. Then there are constants $A,B,C$ that only depend on the fixed data, so that on the fiber $X_y$ and any $0<t\leq 1$ we have
\begin{equation}\label{c2}
\frac{t}{Ce^{Ae^{B\sigma(y)^{-\lambda}}}}\omega_y\leq\ti{\omega}_y\leq tCe^{Ae^{B\sigma(y)^{-\lambda}}}\omega_y,
\end{equation}
\begin{equation}\label{c3}
|\nabla \ti{\omega}_y|^2_{\omega_y}\leq t^{1/2}Ce^{Ae^{B\sigma(y)^{-\lambda}}},
\end{equation}
where $\nabla$ is the covariant derivative of $\omega_y$.
In particular the metrics $\ti{\omega}_y$ converge to zero in $C^1(\omega_y)$ as $t$ approaches zero, uniformly as $y$ varies in a compact set of $Y\backslash f(S)$.
\end{thm}

\begin{thm}\label{weilp2}
As $t\to 0$ the Ricci-flat metrics $\ti{\omega}_t$ on $X\backslash S$ converge to $f^*\omega$, where $\omega$ is a smooth K\"ahler metric on $Y\backslash f(S)$. The convergence is weakly as currents and also in the $C^{1,\beta}_{loc}$ topology of K\"ahler potentials for any $0<\beta<1$. The metric $\omega$ satisfies
$$\Ric(\omega)=\omega_{WP},$$
on $Y\backslash f(S)$, where $\omega_{WP}$ is the pullback of the Weil-Petersson metric from the moduli space of the Calabi-Yau fibers, and it measures the change of complex structures of the fibers.
\end{thm}

\section{A priori estimates}\label{sect3}
In this section we prove {\em a priori} $C^2$ estimates for the degenerating complex Monge-Amp\`ere equations that we are considering, and we also prove $C^3$ estimates along the fibers of $f$.\\

We start with a few lemmas.

\begin{lemma}
There is a uniform constant $C$ so that for all $0<t\leq 1$ we have
\begin{equation}\label{schw}
\tr{\ti{\omega}_t}{\omega_0}\leq C.
\end{equation}
\end{lemma}
\begin{proof}
Recall that we are assuming that $\omega_0= f^*\omega_Z$ where
$f:X\to Z$ is a holomorphic map. We can then use the Chern-Lu formula that appears in Yau's Schwarz lemma computation \cite{yauschw, schwarz} and get
$$\Delta_{\ti{\omega}_t}\log\tr{\ti{\omega}_t}{\omega_0}\geq -A\tr{\ti{\omega}_t}{\omega_0},$$
for a uniform constant $A$. Noticing that
$$\Delta_{\ti{\omega}_t}\varphi_t=n-\tr{\ti{\omega}_t}{\omega_t}\leq 
n-\tr{\ti{\omega}_t}{\omega_0},$$
we see that
\begin{equation}\label{comp1}
\Delta_{\ti{\omega}_t}(\log\tr{\ti{\omega}_t}{\omega_0}-(A+1)\varphi_t)\geq
\tr{\ti{\omega}_t}{\omega_0}-n(A+1).
\end{equation}
Then the maximum principle applied to \eqref{comp1}, together with the estimate \eqref{linftybound}, gives \eqref{schw}.
\end{proof}

The next lemma, which gives a Sobolev constant bound, is due independently to Allard \cite{allard} and Michael-Simon \cite{MS}.
\begin{lemma}\label{sobolevl}
There is a uniform constant $C$ so that for any $0<t\leq 1$, for any $y\in Y\backslash f(S)$ and for any $u\in C^\infty(X_y)$ we have
\begin{equation}\label{sobolev}
\left(\int_{X_y} |u|^{\frac{2(n-m)}{n-m-1}}\omega_y^{n-m}\right)^{\frac{n-m-1}{n-m}}
\leq C \int_{X_y}(|\nabla u|^2_{\omega_y}+|u|^2)\omega_y^{n-m}.
\end{equation}
\end{lemma}
\begin{proof}
For any $y\in Y\backslash f(S)$ the fiber $X_y$ is a smooth $(n-m)$ - dimensional complex submanifold of $X$. Since $X$ is K\"ahler, it follows that $X_y$ is a minimal submanifold, and so it has vanishing mean curvature vector. We then use
the Nash embedding theorem to isometrically embed $(X,\omega_X)$ into Euclidean space, and so we have an isometric embedding $X\to \R^{N}$. The length of the mean curvature vector of the composite isometric embedding $X_y\to X\to \R^{N}$ is then uniformly bounded independent of $y$, since it depends only on the second fundamental form of $X\to \R^N$.
 Then \eqref{sobolev} follows from the uniform Sobolev inequality of \cite{allard, MS}.
Notice that they prove an $L^1$ Sobolev inequality, but this implies the stated $L^2$ Sobolev inequality thanks to the H\"older inequality.
\end{proof}

One can easily avoid the Nash embedding theorem by using a partition of unity to reduce directly to the Euclidean case, but the above proof is perhaps cleaner. 

We note here that the volume of  $X_y$ with respect to $\omega_y$, $\int_{X_y}\omega_y^{n-m}$, is a homological constant independent of $y\in Y\backslash f(S)$, and up to scaling $\omega_X$ we may assume that it is equal to $1$. The next step is to prove a diameter bound for $\omega_y$:

\begin{lemma}\label{fiberdiam}
There is a uniform constant $C$ so that for any $0<t\leq 1$, for any $y\in Y\backslash f(S)$ we have
\begin{equation}\label{fibdiam}
\diam(X_y,\omega_y)\leq C.
\end{equation}
\end{lemma}
\begin{proof}
As above we embed $(X,\omega_X)$ isometrically into $\R^N$ and we get that the length of the mean curvature vector of the composite isometric embedding $X_y\to X\to \R^{N}$ is then uniformly bounded independent of $y$. We can then apply Theorem 1.1 of \cite{topping} and get the required diameter bound. 

Alternatively, first one observes that \eqref{sobolev} implies that there is a uniform constant $\kappa$ so that that geodesic balls in $X_y$ of radius $r<1$ have volume at least $\kappa r^{2(n-m)}$ (Lemma 3.2 in \cite{hebey}). Since the total volume of $X_y$ is constant equal to $1$, an elementary argument gives the required diameter bound.
\end{proof}

The next step is to prove a Poincar\'e inequality for the restricted metric $\omega_y$. This time the constant will not be uniformly bounded, but it will blow up like a power of $\frac{1}{H}$, where $H$ is defined in \eqref{jacob}. To this end, we first estimate the Ricci curvature of $\omega_y$. Fix a point $y\in Y\backslash f(S)$ and choose local coordinates $z^1,\dots,z^{n-m}$ on the fiber $X_y$, which extend locally to coordinates in a ball in $X$. Then pick local coordinates $w^{n-m+1},\dots,w^{n}$ near $y\in Y\backslash f(S)$, so that $z^1,\dots,z^{n-m},z^{n-m+1}=f^*(w^{n-m+1}),\dots,z^n=f^*(w^n)$ give local holomorphic coordinates on $X$. We can also assume that at the point $y$ the metric $\omega_Y$ is the identity.
At any fixed point of $X_y$ we then have
\begin{equation}\label{lowerr}\begin{split}
\Ric(\omega_y)&=-\mn\de\db\log \frac{\omega_y^{n-m}}{dz^1\wedge\dots\wedge
d\ov{z}^{n-m}}\\
&=-\mn\de\db\log \frac{\omega_X^{n-m}\wedge\omega_0^m}{dz^1\wedge\dots\wedge
d\ov{z}^{n}}\\
& =-\mn\de\db\log H-\mn\de\db\log\frac{\omega_X^n}{dz^1\wedge\dots\wedge
d\ov{z}^{n}}\\
&\geq -\frac{\mn\de\db H}{H}+\Ric(\omega_X)|_{X_y}\\
&\geq -\left(\frac{C}{H}+C\right)\omega_y\geq -\frac{C}{H}\omega_y,
\end{split}
\end{equation}
where all derivatives are in fiber directions.
Combining \eqref{lowerr} and \eqref{hack} we see that the Ricci curvature of $\omega_y$ is bounded below by $-C\sigma^{-\lambda}$. Since the diameter of $\omega_y$ is bounded by Lemma \ref{fiberdiam}, a theorem of Li-Yau \cite{LY} then shows that the Poincar\'e constant of $\omega_y$ is bounded above by
$Ce^{B\sigma^{-\lambda}}$. This proves the following

\begin{lemma}\label{poincarel}
There are uniform constants $\lambda, B,C$ so that for any $0<t\leq 1$, for any $y\in Y\backslash f(S)$ and for any $u\in C^\infty(X_y)$ with $\int_{X_y}u\omega_y^{n-m}=0$ we have
\begin{equation}\label{poincare}
\int_{X_y} |u|^{2}\omega_y^{n-m}
\leq Ce^{B\sigma^{-\lambda}} \int_{X_y}|\nabla u|^2_{\omega_y}\omega_y^{n-m}.
\end{equation}
\end{lemma}

We now let $\ti{\omega}_y$ be the restriction $\ti{\omega}_t|_{X_y}$. We have the following estimate for the volume form of $\ti{\omega}_y$ on $X_y$:
\begin{equation}\label{volform}
\begin{split}
\frac{\ti{\omega}_y^{n-m}}{\omega_y^{n-m}}&=\frac{\ti{\omega}_t^{n-m}\wedge\omega_0^m}
{\omega_X^{n-m}\wedge\omega_0^m}=\frac{\ti{\omega}_t^{n-m}\wedge\omega_0^m}{\ti{\omega}_t^n}
\cdot\frac{\ti{\omega}_t^n}{H\omega_X^n}\\
&\leq \left(\frac{\ti{\omega}_t^{n-1}\wedge\omega_0}{\ti{\omega}_t^n}\right)^m \frac{c_t t^{n-m}e^E}{H}\\
&=(\tr{\ti{\omega}_t}{\omega_0})^m \frac{c_t t^{n-m}e^E}{H}
\leq \frac{Ct^{n-m}}{\sigma^\lambda}.
\end{split}
\end{equation}
Notice that when we restrict to $X_y$ we have
$$\ti{\omega}_y=(\omega_0+t\omega_X+\mn\de\db\varphi_t)|_{X_y}=t\omega_y+(
\mn\de\db\varphi_t)|_{X_y}.$$
It is convenient to define a function $\underline{\varphi_t}$ on $Y\backslash f(S)$ 
by
$$\underline{\varphi_t}(y)=\int_{X_y}\varphi_t \omega_y^{n-m}.$$
 This is just the ``integration along the fibers'' of $\varphi_t$, and we will also denote by $\underline{\varphi_t}$ its pullback to $X\backslash S$ via $f$. We also define a function on $X\backslash S$ by
$$\psi=\frac{1}{t}\left(\varphi_t-\underline{\varphi_t}\right),$$ so that we have
$\int_{X_y}\psi\omega_y^{n-m}=0$ and on $X_y$ we have
\begin{equation}\label{maineq}
(\omega_y+\mn\de\db\psi)^{n-m}=\frac{\ti{\omega}_y^{n-m}}{t^{n-m}}\leq \frac{C}{\sigma^\lambda}\omega_y^{n-m}. 
\end{equation}
We can then apply Yau's $L^\infty$ estimate for complex Monge-Amp\`ere equations \cite{Ya} to the inequality \eqref{maineq}. Since the volume of $X_y$ is constant equal to $1$, the Sobolev constant of $\omega_y$ is uniformly bounded (Lemma \ref{sobolevl}) and the Poincar\'e constant is controlled by Lemma \ref{poincarel}, Yau's $L^\infty$ estimate gives
\begin{equation}\label{linft}
\sup_{X_y}\left|\varphi_t-\underline{\varphi_t}\right|=t\sup_{X_y}|\psi|\leq
t Ce^{B\sigma(y)^{-\lambda}},
\end{equation}
where we increased the constant $B$ to absorb the term $\sigma^{-\lambda}$ in \eqref{maineq}.
Recall that from \eqref{linftybound} we have a uniform bound for the oscillation of $\varphi_t$. 

\begin{proof}[Proof of Theorem \ref{estimates}]
First we will show the right-hand side inequality in \eqref{c2a}.
We will apply the maximum principle to the quantity $$K=e^{-B\sigma^{-\lambda}}\left(\log\tr{\omega_X}{\ti{\omega}_t}-\frac{A}{t}(\varphi_t-\underline{\varphi_t})\right),$$
where $A$ is a suitably chosen uniform large constant. The maximum of $K$ on $X\backslash S$ is obviously achieved, and we will show that $K\leq C$ for a uniform constant $C$. This together with \eqref{linft} will show that on
$X\backslash S$ we have
\begin{equation}\label{oneway}
\Delta_{\omega_X}\varphi_t=\tr{\omega_X}{\ti{\omega}_t}-\tr{\omega_X}{\omega_0}-nt\leq
\tr{\omega_X}{\ti{\omega}_t}\leq Ce^{Ce^{B\sigma^{-\lambda}}},
\end{equation}
which is half of \eqref{c2a}
To do this, we first compute as in Yau's $C^2$ estimates \cite{Ya}
$$\Delta_{\ti{\omega}_t}\log\tr{\omega_X}{\ti{\omega}_t}\geq -C\tr{\ti{\omega}_t}{\omega_X}-C,$$
for a uniform constant $C$. On the other hand
$$\Delta_{\ti{\omega}_t}\varphi_t\leq n-t\cdot\tr{\ti{\omega}_t}{\omega_X},$$
and so if $A$ is large enough we get
$$\Delta_{\ti{\omega}_t}\left(\log\tr{\omega_X}{\ti{\omega}_t}-\frac{A}{t}\varphi_t\right)\geq
\tr{\ti{\omega}_t}{\omega_X}-\frac{C}{t}.$$
Since $f$ is locally a submersion on $X\backslash S$, the fiber integration formula
$$\de\db\underline{\varphi_t}=f_*(\de\db\varphi_t\wedge\omega_X^{n-m})$$
holds. So we can compute that
\begin{equation}\begin{split}
\Delta_{\ti{\omega}_t}\underline{\varphi_t}&=\tr{\ti{\omega}_t}
{f_*(\mn\de\db\varphi_t\wedge\omega_X^{n-m})}\\
&=\tr{\ti{\omega}_t}
{f_*((\ti{\omega}_t-\omega_t)\wedge\omega_X^{n-m})}\\
&\geq-\tr{\ti{\omega}_t}
{f_*(\omega_t\wedge\omega_X^{n-m})}\\
&=-\tr{\ti{\omega}_t}
{f_*(f^*\omega_Y\wedge\omega_X^{n-m})}-t\tr{\ti{\omega}_t}
{f_*(\omega_X^{n-m+1})}\\
&=-\tr{\ti{\omega}_t}{\omega_0}-t\tr{\ti{\omega}_t}
{f_*(\omega_X^{n-m+1})}.
\end{split}
\end{equation}
On $Y\backslash f(S)$ the K\"ahler form $f_*(\omega_X^{n-m+1})$ can be estimated by
\begin{equation}\label{utile}
\begin{split}
f_*(\omega_X^{n-m+1})&\leq\frac{\omega_Y^{m-1}\wedge f_*(\omega_X^{n-m+1})}{\omega_Y^m}\omega_Y=\frac{f_*(\omega_0^{m-1}\wedge\omega_X^{n-m+1})}{\omega_Y^m}
\omega_Y\\
&\leq C\frac{f_*(\omega_X^{n})}{\omega_Y^m}\omega_Y=C\frac{f_*(H^{-1}\omega_0^m\wedge\omega_X^{n-m})}{\omega_Y^m}\omega_Y\\
&\leq C \sigma^{-\lambda}\frac{f_*(\omega_0^m\wedge\omega_X^{n-m})}{\omega_Y^m}\omega_Y
=C \sigma^{-\lambda} \omega_Y.
\end{split}
\end{equation}
and so using \eqref{schw} we get
$$\Delta_{\ti{\omega}_t}\underline{\varphi_t}\geq -C-tC\sigma^{-\lambda}.$$ It follows that
\begin{equation}\label{stima1}
\Delta_{\ti{\omega}_t}\left(\log\tr{\omega_X}{\ti{\omega}_t}-\frac{A}{t}(\varphi_t-\underline{\varphi_t})\right)\geq
\tr{\ti{\omega}_t}{\omega_X}-\frac{C}{t}-C\sigma^{-\lambda}.
\end{equation}
Using \eqref{hack3} and \eqref{schw} we have that
\begin{equation}\label{calcc1}
|\Delta_{\ti{\omega}_t}\sigma|\leq C\tr{\ti{\omega}_t}{\omega_0}\leq C,
\end{equation}
\begin{equation}\label{calcc2}
|\nabla\sigma|^2_{\ti{\omega}_t}\leq C\tr{\ti{\omega}_t}{\omega_0}\leq C,
\end{equation}
Using \eqref{stima1} we then compute
\begin{equation}\label{stima2}\begin{split}
\Delta_{\ti{\omega}_t}K&\geq e^{-B\sigma^{-\lambda}}\left(\tr{\ti{\omega}_t}{\omega_X}-\frac{C}{t}-C\sigma^{-\lambda}\right)\\
&+\left(\log\tr{\omega_X}{\ti{\omega}_t}-\frac{A}{t}(\varphi_t-\underline{\varphi_t})\right)
\Delta_{\ti{\omega}_t}\left(e^{-B\sigma^{-\lambda}}\right)\\
&+2e^{B\sigma^{-\lambda}}\mathrm{Re}\langle\nabla K,\nabla e^{-B\sigma^{-\lambda}}\rangle_{\ti{\omega}_t}\\
&-2\left(\log\tr{\omega_X}{\ti{\omega}_t}-\frac{A}{t}(\varphi_t-\underline{\varphi_t})\right)
e^{B\sigma^{-\lambda}}|\nabla e^{-B\sigma^{-\lambda}}|^2_{\ti{\omega}_t}.
\end{split}
\end{equation}
Using \eqref{hack3}, \eqref{calcc1} and \eqref{calcc2}, the second term in \eqref{stima2} can be estimated as follows
\begin{equation}\label{stima3}\begin{split}
\Delta_{\ti{\omega}_t}\left(e^{-B\sigma^{-\lambda}}\right)&=
\frac{B\lambda e^{-B\sigma^{-\lambda}}}{\sigma^{\lambda+1}}\Delta_{\ti{\omega}_t}\sigma
+\frac{B^2\lambda^2 e^{-B\sigma^{-\lambda}}}{\sigma^{2\lambda+2}}|\nabla\sigma|^2_{\ti{\omega}_t}\\
&-\frac{B\lambda(\lambda+1) e^{-B\sigma^{-\lambda}}}{\sigma^{\lambda+2}}|\nabla\sigma|^2_{\ti{\omega}_t}\\
&\geq -C\frac{e^{-B\sigma^{-\lambda}}}{\sigma^{\lambda+1}}-C\frac{e^{-B\sigma^{-\lambda}}}{\sigma^{\lambda+2}}\\
&\geq-C\frac{e^{-B\sigma^{-\lambda}}}{\sigma^{\lambda+2}}.
\end{split}
\end{equation}
At the maximum of $K$ we may assume that $K\geq 0$, otherwise we have nothing to prove. Hence we can use \eqref{linft} to estimate
\begin{equation}\label{stima4}\begin{split}
\left(\log\tr{\omega_X}{\ti{\omega}_t}-\frac{A}{t}(\varphi_t-\underline{\varphi_t})\right)
\Delta_{\ti{\omega}_t}\left(e^{-B\sigma^{-\lambda}}\right)\\
\geq-C\frac{e^{-B\sigma^{-\lambda}}}{\sigma^{\lambda+2}}\log\tr{\omega_X}{\ti{\omega}_t}-
\frac{C}{\sigma^{\lambda+2}}.
\end{split}
\end{equation}
The fourth term in \eqref{stima2} can be estimated using \eqref{calcc2}
\begin{equation}\label{stima5}\begin{split}
|\nabla e^{-B\sigma^{-\lambda}}|^2_{\ti{\omega}_t}&=\frac{B^2\lambda^2
e^{-2B\sigma^{-\lambda}}}{\sigma^{2\lambda+2}}|\nabla\sigma|^2_{\ti{\omega}_t}
\leq \frac{C e^{-2B\sigma^{-\lambda}}}{\sigma^{2\lambda+2}},
\end{split}
\end{equation}
\begin{equation}\label{stima6}\begin{split}
-2\left(\log\tr{\omega_X}{\ti{\omega}_t}-\frac{A}{t}(\varphi_t-\underline{\varphi_t})\right)
e^{B\sigma^{-\lambda}}|\nabla e^{-B\sigma^{-\lambda}}|^2_{\ti{\omega}_t}\\
\geq
-C\frac{e^{-B\sigma^{-\lambda}}}{\sigma^{2\lambda+2}}\log\tr{\omega_X}{\ti{\omega}_t}
-\frac{C}{\sigma^{2\lambda+2}}.
\end{split}
\end{equation}
Plugging \eqref{stima4} and \eqref{stima6} in \eqref{stima2}, at the maximum point of $K$ we get
$$0\geq\tr{\ti{\omega}_t}{\omega_X}-\frac{C}{t}-\frac{C}{\sigma^{\lambda}}
-\frac{C}{\sigma^{2\lambda+2}}\log\tr{\omega_X}{\ti{\omega}_t}-
C\frac{e^{B\sigma^{-\lambda}}}{\sigma^{2\lambda+2}}
.$$
Since for any two K\"ahler metrics $\omega, \ti{\omega}$ we have
\begin{equation}\label{hack2}
\tr{\omega}{\ti{\omega}}\leq (\tr{\ti{\omega}}{\omega})^{n-1}\frac{\ti{\omega}^n}{\omega^n},
\end{equation}
we see that
$$\tr{\omega_X}{\ti{\omega}_t}\leq Ct^{n-m} (\tr{\ti{\omega}_t}{\omega_X})^{n-1}\leq C(\tr{\ti{\omega}_t}{\omega_X})^{n-1},$$
and using this and the inequalities $2ab\leq \ve a^2+b^2/\ve$ and $(\log x)^2\leq x+C$ we get
$$\tr{\ti{\omega}_t}{\omega_X}\leq \frac{C}{t}+\frac{C}{\sigma^{\lambda}}+
\frac{C}{\sigma^{4\lambda+4}}+C\frac{e^{B\sigma^{-\lambda}}}{\sigma^{2\lambda+2}}+\frac{1}{2}\tr{\ti{\omega}_t}{\omega_X},$$
whence
$$\tr{\ti{\omega}_t}{\omega_X}\leq \frac{C}{t}+Ce^{C\sigma^{-\lambda}}.$$
At the same point we then get
$$\tr{\ti{\omega}_t}{\omega_t}=\tr{\ti{\omega}_t}{(\omega_0+t\omega_X)}\leq 
C+tCe^{C\sigma^{-\lambda}}.$$
and using \eqref{hack2} we get
\begin{equation}\label{stima7}
\tr{\omega_t}{\ti{\omega}_t}\leq (\tr{\ti{\omega}_t}{\omega_t})^{n-1}\frac{\ti{\omega}^n_t}{\omega^n_t}\leq
\left(C+tCe^{C\sigma^{-\lambda}}\right)^{n-1}\frac{\ti{\omega}^n_t}{\omega^n_t}.
\end{equation}
We now use \eqref{jacob}, \eqref{degma2} and \eqref{hack} to get
\begin{equation}\label{stima8}
\frac{\ti{\omega}^n_t}{\omega^n_t}\leq \frac{Ct^{n-m}\omega_X^n}{\omega_0^m\wedge(t\omega_X)^{n-m}}=\frac{C}{H}\leq\frac{C}{\sigma^\lambda}.
\end{equation}
Combining \eqref{stima7} and \eqref{stima8} we get
$$\tr{\omega_t}{\ti{\omega}_t}\leq Ce^{C\sigma^{-\lambda}},$$
for some uniform constant $C$. But we also have $\omega_t=\omega_0+t\omega_X\leq C\omega_X$ and so we get
$$\tr{\omega_X}{\ti{\omega}_t}\leq Ce^{C\sigma^{-\lambda}}.$$
Using \eqref{linft} again, this implies that at the maximum of $K$ we have
$$K\leq C+e^{-B\sigma^{-\lambda}}\log (Ce^{C\sigma^{-\lambda}})\leq C.$$

We now show the left-hand side inequality in \eqref{c2a}.
To this extent we apply the maximum principle to the quantity $$K_1=e^{-B\sigma^{-\lambda}_h}\left(\log(t\cdot\tr{\ti{\omega}_t}{\omega_X})-\frac{A}{t}(\varphi_t-\underline{\varphi_t})\right),$$
where $A$ is a suitably chosen uniform large constant. The maximum of $K_1$ on $X\backslash S$ is obviously achieved, and we will show that $K_1\leq C$ for a uniform constant $C$. This together with \eqref{linft} will show that on
$X\backslash S$ we have
\begin{equation}\label{theother}
\tr{\ti{\omega}_t}{\omega_X}\leq  \frac{C}{t}e^{Ce^{B\sigma^{-\lambda}}},
\end{equation}
which is the other half of \eqref{c2a}.
To prove that $K_1\leq C$ we use the maximum principle and,
as in \eqref{stima2}, we compute
\begin{equation}\label{stima11}\begin{split}
\Delta_{\ti{\omega}_t}K_1&\geq e^{-B\sigma^{-\lambda}}\left(\tr{\ti{\omega}_t}{\omega_X}-\frac{C}{t}-C\sigma^{-\lambda}\right)\\
&+\left(\log(t\cdot\tr{\ti{\omega}_t}{\omega_X})-\frac{A}{t}(\varphi_t-\underline{\varphi_t})\right)
\Delta_{\ti{\omega}_t}\left(e^{-B\sigma^{-\lambda}}\right)\\
&+2e^{B\sigma^{-\lambda}}\mathrm{Re}\langle\nabla K_1,\nabla e^{-B\sigma^{-\lambda}}\rangle_{\ti{\omega}_t}\\
&-2\left(\log(t\cdot\tr{\ti{\omega}_t}{\omega_X})-\frac{A}{t}(\varphi_t-\underline{\varphi_t})\right)
e^{B\sigma^{-\lambda}}|\nabla e^{-B\sigma^{-\lambda}}|^2_{\ti{\omega}_t}.
\end{split}
\end{equation}
We estimate this in the same way as before and get
\begin{equation}\label{stima12}\begin{split}
\Delta_{\ti{\omega}_t}K_1&\geq 
e^{-B\sigma^{-\lambda}}\left(\tr{\ti{\omega}_t}{\omega_X}-\frac{C}{t}-C\sigma^{-\lambda}\right)\\
&-C\frac{e^{-B\sigma^{-\lambda}}}{\sigma^{2\lambda+2}}\log(t\cdot\tr{\ti{\omega}_t}{\omega_X})
-\frac{C}{\sigma^{2\lambda+2}}\\
&+2e^{B\sigma^{-\lambda}}\mathrm{Re}\langle\nabla K_1,\nabla e^{-B\sigma^{-\lambda}}\rangle_{\ti{\omega}_t}.
\end{split}
\end{equation}
At the maximum of $K_1$ we get
$$0\geq\tr{\ti{\omega}_t}{\omega_X}-\frac{C}{t}-\frac{C}{\sigma^{\lambda}}
-\frac{C}{\sigma^{2\lambda+2}}\log(t\cdot\tr{\ti{\omega}_t}{\omega_X})-Ce^{C\sigma^{-\lambda}},$$
and using the inequalities $2ab\leq \ve a^2+b^2/\ve$ and $(\log x)^2\leq x+C$ we get
$$\tr{\ti{\omega}_t}{\omega_X}\leq \frac{C}{t}+\frac{C}{\sigma^{\lambda}}+
\frac{C}{\sigma^{4\lambda+4}}+Ce^{C\sigma^{-\lambda}}+\frac{1}{2}\tr{\ti{\omega}_t}{\omega_X},$$
whence
$$t\cdot\tr{\ti{\omega}_t}{\omega_X}\leq C+tCe^{C\sigma^{-\lambda}}\leq Ce^{C\sigma^{-\lambda}},$$
and so at that point 
$$K_1\leq C+e^{-B\sigma^{-\lambda}}\log(Ce^{C\sigma^{-\lambda}})\leq C,$$
and we are done.
\end{proof}
\begin{proof}[Proof of Theorem \ref{estimates2}] 
We will first show \eqref{c2}, which is an easy consequence of \eqref{c2a}.
The left-hand side follows immediately from \eqref{c2a}, which implies
\begin{equation}\label{stima9}
\tr{\ti{\omega}_y}{\omega_y}\leq \frac{C}{t}e^{Ce^{B\sigma^{-\lambda}}}.
\end{equation}
Then \eqref{hack2} and \eqref{volform} give
\begin{equation}\label{stima10}
\tr{\omega_y}{\ti{\omega}_y}\leq (\tr{\ti{\omega}_y}{\omega_y})^{n-m-1}\frac{\ti{\omega}_y^{n-m}}{\omega_y^{n-m}}\leq
t \frac{Ce^{Ce^{B\sigma^{-\lambda}}}}{\sigma^\lambda}\leq t Ce^{Ce^{B\sigma^{-\lambda}}},
\end{equation}
which proves \eqref{c2}.

Next, we show \eqref{c3}. Recall from \eqref{oneway} and \eqref{theother} that on $X\backslash S$ we have
\begin{equation}\label{oneway2}
\tr{\omega_X}{\ti{\omega}_t}\leq Ce^{C_0 e^{B\sigma^{-\lambda}}}.
\end{equation}
\begin{equation}\label{theother2}
\tr{\ti{\omega}_t}{\omega_X}\leq  \frac{C}{t}e^{C_0 e^{B\sigma^{-\lambda}}},
\end{equation}
for uniform constants $B, C, C_0$.
We apply the maximum principle to the quantity 
$$K_2=e^{-Ae^{B\sigma^{-\lambda}}}\left(\Ss+C\frac{e^{3C_0 e^{B\sigma^{-\lambda}}}}{t^{5/2}}\tr{\omega_X}{\ti{\omega}_t}\right),$$
for suitable constants $A, C$,
where the quantity $\Ss$ is the same quantity as in \cite{Ya}:
$$\Ss=|\nabla \ti{\omega}_t|^2_{\ti{\omega}_t},$$ where $\nabla$ is the covariant derivative associated to the metric $\omega_X$. Using $\varphi_t$ we can write
$$\Ss=\ti{g}_t^{i\ov{p}}\ti{g}_t^{q\ov{j}}\ti{g}_t^{k\ov{r}}(g^0_{i\ov{j},k}+\varphi_{i\ov{j}k})
(g^0_{q\ov{p},\ov{r}}+\varphi_{\ov{p}q\ov{r}}),$$
where again lower indices are covariant derivatives with respect to $\omega_X$,
and where $g^0_{i\ov{j}}$ are the components of $\omega_0$.
We are going to show that $K_2\leq \frac{C}{t^{5/2}}$, and using \eqref{oneway2} this implies that
\begin{equation}\label{stima13}
\Ss\leq \frac{Ce^{Ae^{B\sigma^{-\lambda}}}}{t^{5/2}}.
\end{equation}
We now use \eqref{stima10}, which says that on $X_y$ we have
\begin{equation}\label{thethird}
\tr{\omega_y}{\ti{\omega}_y}\leq t Ce^{C_0e^{B\sigma^{-\lambda}}},
\end{equation}
At any given point of $X_y$ we can assume that $\omega_X$ is the identity and $\ti{\omega}_t$ is diagonal with positive entries $\lambda_i$, $1\leq i\leq n$,
so that the first $n-m$ directions are tangent to the fiber $X_y$. Then \eqref{thethird} gives that
\begin{equation}\label{thefourth}
\lambda_i\leq  t Ce^{C_0e^{B\sigma^{-\lambda}}},
\end{equation}
for $1\leq i\leq n-m$. Then using \eqref{stima13} we see that
$$\sum_{i,j,k=1}^{n-m}\frac{1}{\lambda_i\lambda_j\lambda_k}
|\varphi_{i\ov{j} k}|^2\leq \sum_{i,j,k=1}^{n}\frac{1}{\lambda_i\lambda_j\lambda_k}
|g^0_{i\ov{j},k}+\varphi_{i\ov{j} k}|^2=\Ss\leq \frac{Ce^{Ae^{B\sigma^{-\lambda}}}}{t^{5/2}},$$
where we have used that $g^0_{i\ov{j},k}$ vanishes whenever $1\leq i,j\leq n-m$
since $\omega_0|_{X_y}=0$. Using \eqref{thefourth} we get
$$|\nabla \ti{\omega}_y|^2_{\omega_y}=\sum_{i,j,k=1}^{n-m} |\varphi_{i\ov{j} k}|^2\leq t^{1/2}Ce^{(A+3C_0)e^{B\sigma^{-\lambda}}},$$
and this is \eqref{c3}.

We now prove that $K_2\leq \frac{C}{t^{5/2}}$. To simplify the computation, we will use the notation
$$\F(x)=e^{xe^{B\sigma^{-\lambda}}},$$
where $x$ is a real number, and we note here that $\F$ is increasing.
The starting point is a precise formula for $\Delta_{\ti{\omega}_t}\Ss$. This is just Yau's $C^3$ estimate \cite{Ya}, but without assuming that the metrics 
$\ti{\omega}_t$ and $\omega_t=\omega_0+t\omega_X$ are equivalent, and it is done 
in a more general setting in \cite{taming} (see also \cite{PSS}). With the notation of \cite{taming} we can write
$$\Ss=\sum_{i,j,k}|a^i_{jk}|^2.$$
We then choose local unitary
frames $\{\theta^1,\dots,\theta^n\}$ for $\omega_X$ and
$\{\ti{\theta}^1,\dots,\ti{\theta}^n\}$ for $\ti{\omega}_t$, and write
$$\ti{\theta}^i=\sum_j a^i_j\theta^j,$$
$$\theta^i=\sum_j b^i_j\ti{\theta}^j,$$
for some local matrices of functions $a^i_j, b^i_j$. Notice that at any given point we can choose the frames and arrange that 
\begin{equation}\label{stima14}
a^i_j=\sqrt{\lambda_i}\delta^i_j,
\end{equation}
\begin{equation}\label{stima15}
b^i_j=\frac{1}{\sqrt{\lambda_i}}\delta^i_j.
\end{equation}
Then in our case \cite[(4.3)]{taming} reads
\begin{eqnarray}\label{stimaa} \nonumber
\Delta_{\ti{\omega}_t}\Ss
 & \geq & 
 2 \mathrm{Re} \biggl( \ov{a^i_{k\ell}} \biggl( b^m_k b^q_\ell \ov{b^s_p} R^j_{mq\ov{s}} a^i_{r p} a^r_j 
- a^i_j b^q_\ell \ov{b^s_p} R^j_{mq\ov{s}} a^r_{kp} b^m_r   \\ 
&&  \mbox{}  - a^i_j b^m_k \ov{b^s_p} R^j_{mq\ov{s}} a^r_{\ell p} b^q_r
 + a^i_j b^m_k b^q_\ell \ov{b^s_p} b_p^u  R^j_{mq\ov{s},u}  \biggr)\biggr), 
 \end{eqnarray}
where we are summing over all indices, $R^j_{mq\ov{s}}$ represents the curvature of $\omega_X$ and $R^j_{mq\ov{s},u}$ its covariant derivative (with respect to $\omega_X$). Since these are fixed tensors, we can use the Cauchy-Schwarz inequality and \eqref{stima14}, \eqref{stima15} to estimate the first term on the right hand side of \eqref{stimaa} by
\begin{equation*}\begin{split}
 &\left|2\mathrm{Re}\left(\ov{a^i_{k\ell}} b^m_k b^q_\ell \ov{b^s_p} R^j_{mq\ov{s}} a^i_{r p} a^r_j \right)\right|\leq C \sum_{i,k,\ell,r,p} 
|a^i_{k\ell} a^i_{r p}|\sqrt{\frac{\lambda_r}{\lambda_k\lambda_\ell\lambda_p}}\\
&\leq C \left(\sum_j\lambda_j\right)^{\frac{1}{2}}\left(\sum_q\frac{1}{\lambda_q}\right)^{\frac{3}{2}}
\sum_{k,\ell,r,p}\left(\sum_i |a^i_{k\ell}|^2\right)^{\frac{1}{2}}\left(\sum_i |a^i_{r p}|^2\right)^{\frac{1}{2}}\\
&= C(\tr{\omega_X}{\ti{\omega}_t})^{\frac{1}{2}}(\tr{\ti{\omega}_t}{\omega_X})^{\frac{3}{2}}
\left(\sum_{i,k,\ell} |a^i_{k\ell}|^2\right)^{\frac{1}{2}}\left(\sum_{i,r,p} |a^i_{r p}|^2\right)^{\frac{1}{2}}\\
&=C \Ss (\tr{\omega_X}{\ti{\omega}_t})^{\frac{1}{2}}(\tr{\ti{\omega}_t}{\omega_X})^{\frac{3}{2}}.
\end{split}\end{equation*}
The second and third term in \eqref{stimaa} are estimated similarly, while the fourth term can be bounded by
\begin{equation*}\begin{split}
&\left|2\mathrm{Re}\left(\ov{a^i_{k\ell}} a^i_j b^m_k b^q_\ell \ov{b^s_p} b_p^u  R^j_{mq\ov{s},u}\right)\right|\leq C \sum_{i,k,\ell,p}|a^i_{k\ell}|
\sqrt{\frac{\lambda_i}{\lambda_k\lambda_\ell\lambda_p^2}}\\
&\leq C\left(\sum_j\lambda_j\right)^{\frac{1}{2}}\left(\sum_q\frac{1}{\lambda_q}\right)^{2}
\sum_{i,k,\ell}|a^i_{k\ell}|\\
&\leq C\sqrt{\Ss}(\tr{\omega_X}{\ti{\omega}_t})^{\frac{1}{2}}(\tr{\ti{\omega}_t}{\omega_X})^{2}.
\end{split}\end{equation*}
Overall we can estimate
\begin{equation}\label{stimab}
\Delta_{\ti{\omega}_t}\Ss\geq -C\Ss (\tr{\ti{\omega}_t}{\omega_X})^{3/2}(\tr{\omega_X}{\ti{\omega}_t})^{1/2}
-C\sqrt{\Ss}(\tr{\ti{\omega}_t}{\omega_X})^{2}(\tr{\omega_X}{\ti{\omega}_t})^{1/2}.
\end{equation}
On the other hand from \cite[Lemma 3.3]{taming} we see that
\begin{equation}\label{stimac}\begin{split}
\Delta_{\ti{\omega}_t}\tr{\omega_X}{\ti{\omega}_t}&=a^i_{k\ell} \ov{a^i_{p\ell}} a^k_j \ov{a^p_j} + \ov{a^i_j} a^i_r b^q_\ell \ov{b^s_\ell} R^r_{jq \ov{s}}\\
&\geq \sum_{i,j,\ell}|a^i_{j\ell}|^2\lambda_j-C\sum_{i,\ell}\frac{\lambda_i}{\lambda_\ell}\\
&\geq \left(\sum_k\frac{1}{\lambda_k}\right)^{-1}\sum_{i,j,\ell}|a^i_{j\ell}|^2
-C\left(\sum_p\lambda_p\right)\left(\sum_q\frac{1}{\lambda_q}\right)\\
&=\frac{\Ss}{\tr{\ti{\omega}_t}{\omega_X}}-C(\tr{\ti{\omega}_t}{\omega_X})(\tr{\omega_X}{\ti{\omega}_t}).
\end{split}\end{equation}
We now insert \eqref{oneway2}, \eqref{theother2} in \eqref{stimab}, \eqref{stimac} and get
\begin{equation}\label{stima18}
\Delta_{\ti{\omega}_t}\Ss\geq -\frac{C\F(2C_0)}{t^{3/2}}\Ss 
-\frac{C\F(5C_0/2)}{t^{2}}\sqrt{\Ss},
\end{equation}
$$\Delta_{\ti{\omega}_t}\tr{\omega_X}{\ti{\omega}_t}\geq \frac{t \F(-C_0)}{C}\Ss-\frac{C\F(2C_0)}{t}.$$
We then compute
\begin{equation}\label{stima16}
\begin{split}
\Delta_{\ti{\omega}_t}\left(\frac{\F(3C_0)}{t^{5/2}}\tr{\omega_X}{\ti{\omega}_t}\right)
&\geq \frac{\F(2C_0)}{Ct^{3/2}}\Ss-\frac{C\F(5C_0)}{t^{7/2}}\\
&+\frac{2}{t^{5/2}}\mathrm{Re}\langle\nabla \F(3C_0),\nabla \tr{\omega_X}{\ti{\omega}_t}\rangle_{\ti{\omega}_t}\\
&+\frac{1}{t^{5/2}}(\tr{\omega_X}{\ti{\omega}_t})\Delta_{\ti{\omega}_t}\F(3C_0),
\end{split}\end{equation}
and estimate 
$$\mathrm{Re}\langle\nabla \F(3C_0),\nabla \tr{\omega_X}{\ti{\omega}_t}\rangle_{\ti{\omega}_t}\geq -|\nabla \F(3C_0)|_{\ti{\omega}_t}
|\nabla\tr{\omega_X}{\ti{\omega}_t}|_{\ti{\omega}_t}.$$
Using \cite[(3.20)]{taming} we see that
$$|\nabla\tr{\omega_X}{\ti{\omega}_t}|_{\ti{\omega}_t}\leq \sqrt{\Ss} (\tr{\omega_X}{\ti{\omega}_t}).$$
On the other hand a direct computation using \eqref{calcc1} and \eqref{calcc2} shows that there is a constant $C$ such that for any real number $x$ we have
$$|\nabla \F(x)|_{\ti{\omega}_t}\leq C \F(x+1),$$
$$|\Delta_{\ti{\omega}_t}\F(x)|\leq C \F(x+1),$$
and so we have
\begin{equation}\label{stima17}
\begin{split}
\Delta_{\ti{\omega}_t}\left(\frac{\F(3C_0)}{t^{5/2}}\tr{\omega_X}{\ti{\omega}_t}\right)
&\geq \frac{\F(2C_0)}{Ct^{3/2}}\Ss-\frac{C\F(5C_0)}{t^{7/2}}
-\frac{C\F(5C_0)}{t^{5/2}}\sqrt{\Ss}
\\
&-\frac{C\F(5C_0)}{t^{5/2}}.
\end{split}\end{equation}
This and \eqref{stima18} give
\begin{equation*}
\begin{split}
\Delta_{\ti{\omega}_t}\left(\Ss+\frac{C\F(3C_0)}{t^{5/2}}\tr{\omega_X}{\ti{\omega}_t}\right)
&\geq \frac{\F(2C_0)}{t^{3/2}}\Ss-\frac{C\F(5C_0/2)}{t^{2}}\sqrt{\Ss}\\
&-\frac{C\F(5C_0)}{t^{7/2}}
-\frac{C\F(5C_0)}{t^{5/2}}\sqrt{\Ss}
-\frac{C\F(5C_0)}{t^{5/2}}\\
&\geq \frac{\F(2C_0)}{t^{3/2}}\Ss-\frac{C\F(5C_0)}{t^{7/2}}
-\frac{C\F(5C_0)}{t^{5/2}}\sqrt{\Ss},
\end{split}\end{equation*}
and
\begin{equation}\label{stima20}
\begin{split}
\Delta_{\ti{\omega}_t}K_2
&\geq \F(-A)\biggl(\frac{\F(2C_0)}{t^{3/2}}\Ss-\frac{C\F(5C_0)}{t^{7/2}}
-\frac{C\F(5C_0)}{t^{5/2}}\sqrt{\Ss}\\
& -C\F(1)\Ss-\frac{C\F(4C_0+1)}{t^{5/2}}\biggr)
+2\F(A)\mathrm{Re}\langle \nabla K_2,\nabla \F(-A)\rangle_{\ti{\omega}_t}\\
&\geq \F(-A)\biggl(\frac{\F(2C_0)}{Ct^{3/2}}\Ss-\frac{C\F(5C_0)}{t^{7/2}}
-\frac{C\F(5C_0)}{t^{5/2}}\sqrt{\Ss}\biggr)\\
&+2\F(A)\mathrm{Re}\langle \nabla K_2,\nabla \F(-A)\rangle_{\ti{\omega}_t}.
\end{split}\end{equation}
At the maximum of $K_2$ we then get
$$\Ss\leq \frac{C\F(3C_0)}{t}\sqrt{\Ss}+\frac{C\F(3C_0)}{t^2},$$
which implies that
$$\Ss\leq \frac{C\F(6C_0)}{t^2},$$
and so
$$K_2=\F(-A)\left(\Ss+\frac{C\F(3C_0)}{t^{5/2}}\tr{\omega_X}{\ti{\omega}_t}\right)\leq \F(-A)\frac{C\F(6C_0)}{t^{5/2}}\leq \frac{C}{t^{5/2}},$$
if we choose $A\geq 6C_0$.
\end{proof}

{\bf Remark.} In the estimates proved in this section we have repeatedly used the fact that the metrics $\ti{\omega}_t$ are Ricci-flat. If instead one is dealing with the general equation \eqref{degma2}, the only estimate that does not generalize immediately is \eqref{schw} (which requires that the Ricci curvature of $\ti{\omega}_t$ be nonnegative). \\

{\bf Remark.} In the context of collapsing of the K\"ahler-Ricci flow on projective manifolds with semi-ample
canonical bundle and positive Kodaira dimension, parabolic analogues of the $C^2$ estimates \eqref{c2a}, \eqref{c2} were proved by Song-Tian in \cite{SoT3}.

\section{Collapsing of Ricci-flat metrics}\label{sect4}
In this section we use the estimates from section \ref{sect3} to prove that
the Ricci-flat metrics collapse to the base of the fibration.\\

We first explain the meaning of the Weil-Petersson metric, following the discussion in \cite{SoT2}. Recall that the Ricci-flat K\"ahler metric on $X$ cohomologous to $\omega_1=\omega_0+\omega_X$ is denoted by $\ti{\omega}_1$. We will call $\Omega=\ti{\omega}_1^n$ its volume form.  The generic fiber $X_y$ of $f$ is an $(n-m)$-dimensional Calabi-Yau manifold,
and it is equipped with the K\"ahler form $\omega_y=\omega_X|_{X_y}$.

Recall that the volume of $X_y$ is a homological constant independent of $y$, and that we assume that it is equal to $1$. Since $c_1(X_y)=0$, there is a smooth function $F_y$ such that $\Ric(\omega_y)=\mn\de\db F_y$ and $\int_{X_y}(e^{F_y}-1)\omega_y^{n-m}=0$. The functions $F_y$ vary smoothly in $y$, since so do the K\"ahler forms $\omega_y$.
By Yau's theorem there is a unique Ricci-flat K\"ahler metric $\omega_{SF,y}$ on $X_y$ cohomologous to $\omega_y$, given by the solution of
\begin{equation}\label{sf}
\omega_{SF,y}^{n-m}=e^{F_y}\omega_y^{n-m}.
\end{equation}
If we write $\omega_{SF,y}=\omega_y+\mn\de\db\zeta_y$, the functions $\zeta_y$
 vary smoothly in $y$ and so they define a smooth function $\zeta$ on $X\backslash S$. We then define a real closed $(1,1)$-form $\omega_{SF}$ on $X\backslash S$ by $\omega_{SF}=\omega_X+\mn\de\db\zeta$, and call it the semi-flat form. Notice that $\omega_{SF}$ is not necessarily nonnegative (it is K\"ahler only in the fiber directions), but on $X\backslash S$ the $(n,n)$-form $\omega_{SF}^{n-m}\wedge\omega_0^m$ is strictly positive, and so we can define a smooth positive function $F$ on $X\backslash S$ by
\begin{equation}\label{eq10}
F=\frac{\Omega}{\omega_{SF}^{n-m}\wedge\omega_0^m}.
\end{equation}
We claim that $F$ is actually constant on each fiber $X_y$, and so it is the pullback of a function on $Y\backslash f(S)$. To see this, fix a point $y\in Y\backslash f(S)$ and choose local coordinates $z^1,\dots,z^{n-m}$ on the fiber $X_y$, which extend locally to coordinates in a ball in $X$. Then take local coordinates $w^{n-m+1},\dots,w^{n}$ near $y\in Y\backslash f(S)$, so that $z^1,\dots,z^{n-m},z^{n-m+1}=f^*(w^{n-m+1}),\dots,z^n=f^*(w^n)$ give local holomorphic coordinates on $X$. In these coordinates write 
$$\omega_0=\mn\sum_{i,j=n-m+1}^{n}g^0_{i\ov{j}}dz^i\wedge d\ov{z}^j,$$
$$\omega_{SF,y}=\mn\sum_{i,j=1}^{n-m}g^{SF}_{i\ov{j}}dz^i\wedge d\ov{z}^j,$$
$$\Omega=G (\mn)^n dz^1\wedge\dots\wedge d\ov{z}^n.$$
Then locally 
$$F=\frac{G}{\det(g^0_{i\ov{j}})\det(g^{SF}_{i\ov{j}})},$$
and so on the fiber $X_y$ we have
$$\mn\de\db\log F=-\Ric(\ti{\omega}_1)+\Ric(\omega_{SF,y})=0,$$
because $\omega_0$ is the pullback of a metric from $Y$, and so $F$ is indeed constant on $X_y$. Moreover, it is easy to check \cite[Lemma 3.3]{SoT2} that on $Y\backslash f(S)$
we have
$$F=\frac{f_*\Omega}{\omega_Y^m},$$
and so 
$$\int_Y F\omega_Y^m=\int_X\Omega=\int_X\omega_1^n
$$ is finite. In fact there is a positive $\ve$ so that
$\int_Y F^{1+\ve}\omega_Y^m$ is finite \cite[Proposition 3.2]{SoT2}. Then we apply \cite[Theorem 3.2]{SoT2}, which relies on the seminal work of Ko\l odziej \cite{kol} and further generalizations \cite{egz1, zhangthesis}, to solve (uniquely) the complex Monge-Amp\`ere equation
\begin{equation}\label{wp}
(\omega_Y+\mn\de\db\psi)^{m}=\frac{\int_X\omega_0^m\wedge\omega_X^{n-m}}{\int_X\omega_1^n}F\omega_Y^m,
\end{equation}
with $\psi\in L^\infty(Y)$
and moreover $\psi$ is smooth on $Y\backslash f(S)$ (the proof of this follows the arguments of Yau in \cite{Ya}).
We will call $\omega=\omega_Y+\mn\de\db\psi$ the K\"ahler metric on $Y\backslash f(S)$ that we've just constructed. Its Ricci curvature is the Weil-Petersson metric that we are about to define. Recall that the fibers $X_y$ have torsion canonical bundle, so that there is a number $k$ such that $K_{X_y}^{\otimes k}$ is trivial for all $y\in Y\backslash f(S)$.
The Weil-Petersson metric is a smooth nonnegative $(1,1)$-form on $Y\backslash f(S)$ defined as the curvature form of a pseudonorm on the relative canonical line bundle $f_* (\Omega^{n-m}_{X/Y})^{\otimes k}$: if $\Psi_y$ is a local nonzero holomorphic section of $f_* (\Omega^{n-m}_{X/Y})^{\otimes k}$, which means that $\Psi_y$ is a nonzero holomorphic $k$-pluricanonical form on $X_y$ that varies holomorphically in $y$, then we let its length be
$$|\Psi_y|^2_{h_{WP}}=\int_{X_y}(\Psi_y\wedge\ov{\Psi_y})^{\frac{1}{k}}.$$
For $k>1$ this is not a Hermitian metric, but just a pseudonorm.
The Weil-Petersson metric $\omega_{WP}$ on $Y\backslash f(S)$ is just formally the curvature of $h_{WP}$, that is locally we set
$$\omega_{WP}=-\mn\de\db\log|\Psi_y|^2_{h_{WP}},$$
and this is well-defined because the bundle $K_{X_y}^{\otimes k}$ is trivial.
It is a classical fact (see \cite{fs}) that $\omega_{WP}$ is pointwise nonnegative. As an aside, we note here that one can realize $\omega_{WP}$ as the actual curvature form
of an honest Hermitian metric on a relative canonical bundle if one takes a finite unramified $k$-sheeted cyclic cover $\ti{X}\to X$ so that the smooth fibers of $\ti{X}\to Y$ now have trivial canonical bundle.

\begin{prop}[cfr. \cite{SoT2}] On $Y\backslash f(S)$ we have
\begin{equation}\label{wp2}
\Ric(\omega)=\omega_{WP}.
\end{equation}
\end{prop}
\begin{proof}
Differentiating \eqref{wp} we see that
$$\Ric(\omega)=\Ric(\omega_Y)-\mn\de\db\log F.$$
If we fix $y\in Y\backslash f(S)$ and choose $\Psi$ a local never vanishing holomorphic section of $f_* (\Omega^{n-m}_{X/Y})^{\otimes k}$, then we can define a local function $u=\frac{(\Psi\wedge\ov{\Psi})^{1/k}}{\omega_{SF}^{n-m}}$ on $X\backslash S$, which is constant on each fiber $X_y$. Since $\int_{X_y}\omega_{SF}^{n-m}=1$, we see that
$$-\mn\de\db \log u = \omega_{WP}.$$
Then
\begin{equation}\label{wp3}
\Ric(\omega)=\Ric(\omega_Y)-\mn\de\db\log \frac{u \Omega}{(\Psi\wedge\ov{\Psi})^{\frac{1}{k}}\wedge\omega_0^m}.
\end{equation}
Picking local coordinates $z^i$ as above, and writing 
$$\Psi=K [(\mn)^{n-m} dz^1\wedge \dots\wedge dz^{n-m}]^{\otimes k},$$
$$\omega_0=\mn\sum_{i,j=n-m+1}^{n}g^0_{i\ov{j}}dz^i\wedge d\ov{z}^j,$$
$$\Omega=G (\mn)^n dz^1\wedge\dots\wedge d\ov{z}^n,$$
we see that
$$\frac{u \Omega}{(\Psi\wedge\ov{\Psi})^{\frac{1}{k}}\wedge\omega_0^m}=\frac{uG}{|K|^{\frac{2}{k}} \det(g^0_{i\ov{j}})},$$
and since $K$ is holomorphic and $\Omega$ is Ricci-flat we see that
$$-\mn\de\db\log\frac{uG}{|K|^{\frac{2}{k}} \det(g^0_{i\ov{j}})}=\omega_{WP}-\Ric(\omega_Y),$$
which together with \eqref{wp3} gives \eqref{wp2}.
\end{proof}
With these preparations, we can now show Theorem \ref{weilp2}, which can be recast as follows
\begin{thm}\label{weilp}
Consider the Ricci-flat metrics $\ti{\omega}_t$ on $X$, which can be written as
$\ti{\omega}_t=\omega_0+t\omega_X+\mn\de\db\varphi_t$. As $t\to 0$, we have $\varphi_t\to\psi$ in the $C^{1,\beta}_{loc}$ topology on $X\backslash S$, for any $0<\beta<1$,
and so $\ti{\omega}_t$ converges in this topology to $\omega$, which satisfies \eqref{wp2}. Moreover $\ti{\omega}_t$ also converge to $\omega$ weakly as currents on $X$.
\end{thm}
\begin{proof}
We first prove that $\ti{\omega}_t$ converges to $\omega$ in the weak topology of currents. Since the cohomology class of $\ti{\omega}_t$ is bounded, weak compactness of currents implies that from any sequence $t_i\to 0$ we can extract a subsequence so that $\ti{\omega}_{t_i}$ converges weakly to a limit closed positive $(1,1)$-current $\hat{\omega}$, which a priori depends on the sequence.
If we write $\hat{\omega}=\omega_0+\mn\de\db\hat{\varphi}$, it follows that 
$\varphi_{t_i}\to \hat{\varphi}$ in $L^1$, and from the bound \eqref{linftybound} we infer that $\hat{\varphi}$ is in $L^\infty$. Moreover
restricting $\hat{\omega}$ to any smooth fiber $X_y$ we see that
$$\mn\de\db\hat{\varphi}|_{X_y}\geq 0,$$
and the maximum principle implies that $\hat{\varphi}$ is constant on each fiber, and so descends to a bounded function $\hat{\varphi}$ on $Y\backslash f(S)$. We will show that $\hat{\varphi}$ satisfies the same equation \eqref{wp}
as $\psi$, and so by uniqueness $\hat{\varphi}=\psi$. To this end we fix an arbitrary compact set $K\subset Y\backslash f(S)$, and we wish to show that $\hat{\varphi}$ satisfies \eqref{wp} on $K$. 

We then fix $\eta$ a smooth function with support contained in $K$, and we will also denote by $\eta$ its pullback to $X$ via $f$. Recall that we have called $\ti{\omega}_1$ the Ricci-flat metric in the class $[\omega_1]$, and $\Omega=\ti{\omega}_1^n$. Then from the Monge-Amp\`ere equation \eqref{degma1} we have
\begin{equation}\label{to5}
\int_X\eta \Omega=\frac{1}{a_t}\int_X\eta (\omega_0+t\omega_X+\mn\de\db\varphi_t)^n,
\end{equation}
where the constants $a_t$ are equal to
$$\frac{\int_X \omega_t^n}{\int_X \omega_1^n},$$
and behave like \eqref{normconst}. 
We can also write
\begin{equation}\label{to6}
\int_X\eta \Omega=\int_X\eta F \omega_{SF}^{n-m}\wedge\omega_0^m.
\end{equation}
We are now going to estimate $\frac{1}{a_t}\int_X\eta (\omega_0+t\omega_X+\mn\de\db\varphi_t)^n.$ We have
\begin{equation*}\begin{split}
&\frac{1}{a_t}\int_X\eta (\omega_0+t\omega_X+\mn\de\db\varphi_t)^n\\
&=
\frac{1}{a_t}\int_X\eta \left((\omega_0+\mn\de\db\underline{\varphi_t})+(t\omega_X+\mn\de\db(\varphi_t-\underline{\varphi_t})\right)^n\\
&=\frac{1}{a_t}\int_X\eta \sum_{k=0}^{n}\binom{n}{k}(\omega_0+\mn\de\db\underline{\varphi_t})^k\wedge(t\omega_X+\mn\de\db(\varphi_t-\underline{\varphi_t}))^{n-k}
\end{split}
\end{equation*}
First of all observe that the form $\omega_0+\mn\de\db\underline{\varphi_t}$ is the pullback of a form on $Y$, and it can be wedged with itself at most $m$ times, so all terms in the sum with $k>m$ are zero. Next, we claim that all the terms with $k<m$ go to zero as $t\to 0$. To see this, start by observing that on the compact set $K$ the estimate \eqref{oneway2} gives a constant $C$ (that depends on $K$) such that
\begin{equation}\label{to1}
-C\omega_X\leq \mn\de\db\varphi_t\leq C\omega_X.
\end{equation}
Moreover from the equation
$$\de\db\underline{\varphi_t}=f_*(\de\db\varphi_t\wedge\omega_X^{n-m})$$
together with \eqref{to1}, \eqref{utile}, we see that on $f(K)$ we have
\begin{equation}\label{to2}
-C\omega_Y\leq \mn\de\db\underline{\varphi_t}\leq C\omega_Y.
\end{equation}
We also need to use \eqref{linft} which on $K$ gives
\begin{equation}\label{to8}
\sup_K |\varphi_t-\underline{\varphi_t}|\leq Ct.
\end{equation}
Then any term with $k<m$ is equal to
$$\frac{\binom{n}{k}}{a_t}\int_X \eta(\omega_0+\mn\de\db\underline{\varphi_t})^k\wedge(t\omega_X+\mn\de\db(\varphi_t-\underline{\varphi_t}))^{n-k},$$
and it can be expanded into
$$\frac{\binom{n}{k}}{a_t}\sum_{i=0}^{n-k}\binom{n-k}{i}\int_X \eta(\omega_0+\mn\de\db\underline{\varphi_t})^k\wedge(t\omega_X)^{n-k-i}\wedge(\mn\de\db(\varphi_t-\underline{\varphi_t}))^{i}.$$
On $K$ the $(1,1)$-form $\omega_0+\mn\de\db\underline{\varphi_t}$ is bounded by \eqref{to2}. Since $a_t=O(t^{n-m})$ from \eqref{normconst}, we see that the term in this sum with $i=0$ goes to zero. Any term with $i>0$ is comparable to
\begin{equation}\label{to0}
\begin{split}
\frac{1}{t^{n-m}}\int_X 
(\varphi_t-\underline{\varphi_t}) 
\mn\de\db\eta\wedge(\omega_0+\mn\de\db\underline{\varphi_t})^k\wedge(t\omega_X)^{n-k-i}\wedge\\
\wedge(\mn\de\db(\varphi_t-\underline{\varphi_t}))^{i-1}.
\end{split}
\end{equation}
Notice that all the $(1,1)$-forms appearing inside the integral are bounded by \eqref{to1}, \eqref{to2}, and that the function $\varphi_t-\underline{\varphi_t}$ is $O(t)$ by \eqref{to8}. On $K$ the estimate \eqref{c2} gives
\begin{equation}\label{to9}
-Ct\omega_y\leq (\mn\de\db\varphi_t)|_{X_y}=(\mn\de\db(\varphi_t-\underline{\varphi_t}))|_{X_y}\leq Ct\omega_y. 
\end{equation}
The form $\mn\de\db\eta\wedge(\omega_0+\mn\de\db\underline{\varphi_t})^k$ is the pullback of a form from $Y$, and so we can use \eqref{to9} to estimate
\begin{equation*}\begin{split}
\left|\frac{\mn\de\db\eta\wedge(\omega_0+\mn\de\db\underline{\varphi_t})^k\wedge(t\omega_X)^{n-k-i}\wedge(\mn\de\db(\varphi_t-\underline{\varphi_t}))^{i-1}}
{\omega_X^n} \right|\leq \\
\leq C t^{n-m},
\end{split}\end{equation*}
and so the term \eqref{to0} goes to zero. This proves our claim. 

We are then left with only the term with $k=m$, which is
$$\frac{1}{a_t}\int_X\eta\binom{n}{m}(\omega_0+\mn\de\db\underline{\varphi_t})^m\wedge(t\omega_X+\mn\de\db(\varphi_t-\underline{\varphi_t}))^{n-m},$$
and if we expand the term $(t\omega_X+\mn\de\db(\varphi_t-\underline{\varphi_t}))^{n-m}$, we get
\begin{equation*}\begin{split}
\frac{1}{a_t}&\int_X\eta\binom{n}{m}(\omega_0+\mn\de\db\underline{\varphi_t})^m\wedge (t\omega_X)^{n-m}\\
&+\frac{1}{a_t}\int_X \mn\de\db\eta\wedge(\omega_0+\mn\de\db\underline{\varphi_t})^m\wedge\dots,
\end{split}\end{equation*}
and the second term is zero because $\de\db\eta$ is the pullback of a form from the base. We are then left with the term
\begin{equation}\label{to3}
\frac{1}{a_t}\int_X\eta\binom{n}{m}(\omega_0+\mn\de\db\underline{\varphi_t})^m\wedge (t\omega_X)^{n-m},
\end{equation}
which we need to further estimate. Using \eqref{to1} we see that, up to taking a further subsequence, the functions
$\varphi_{t_i}$ converge to $\hat{\varphi}$ in the $C^{1,\beta}(K)$ topology,
and \eqref{to8} implies that the functions $\underline{\varphi_{t_i}}$ also converge to 
$\hat{\varphi}$ uniformly. 
We can then rewrite \eqref{to3} as
$$\frac{t^{n-m}\binom{n}{m}}{a_t}\int_X\eta(\omega_0+\mn\de\db\underline{\varphi_t})^m\wedge \omega_X^{n-m}.$$
Using \eqref{normconst} we see that as $t$ goes to zero the coefficient $\frac{t^{n-m}\binom{n}{m}}{a_t}$ converges to
$$\frac{\int_X\omega_1^n}{\int_X\omega_0^m\wedge\omega_X^{n-m}}.$$
On the other hand we have
\begin{equation*}\begin{split}
\int_X\eta(\omega_0+\mn\de\db\underline{\varphi_t})^m&\wedge \omega_X^{n-m}\\
&=\sum_{k=0}^m\binom{m}{k}\int_X\eta \omega_0^{m-k}\wedge (\mn\de\db\underline{\varphi_t})^k\wedge \omega_X^{n-m}.
\end{split}\end{equation*}
The term with $k=0$ is independent of $t$, while any term with $k>0$ can be written as
\begin{equation}\label{to4}
\int_X \underline{\varphi_t} \mn\de\db\eta\wedge \omega_0^{m-k}\wedge (\mn\de\db\underline{\varphi_t})^{k-1}\wedge \omega_X^{n-m}.
\end{equation}
The $(n,n)$-form $\mn\de\db\eta\wedge \omega_0^{m-k}\wedge (\mn\de\db\underline{\varphi_t})^{k-1}\wedge \omega_X^{n-m}$ is supported in $K$ and is uniformly bounded by \eqref{to2}, and the functions $\underline{\varphi_{t_i}}$ converge uniformly to $\hat{\varphi}$, and so along the sequence $t_i$ the term \eqref{to4} has the same limit as
$$\int_X \hat{\varphi} \mn\de\db\eta\wedge \omega_0^{m-k}\wedge (\mn\de\db\underline{\varphi_t})^{k-1}\wedge \omega_X^{n-m}.$$
But this is equal to
$$\int_X
\underline{\varphi_t} \mn\de\db\eta\wedge \omega_0^{m-k}\wedge (\mn\de\db\underline{\varphi_t})^{k-2}\wedge\mn\de\db\hat{\varphi} \wedge \omega_X^{n-m},$$
and repeating the same argument $k-1$ times we see that along the sequence $t_i$
the term \eqref{to4} converges to
$$\int_X\eta \omega_0^{m-k}\wedge (\mn\de\db\hat{\varphi})^k\wedge \omega_X^{n-m}.$$
It follows that along the sequence $t_i$ the term \eqref{to3} converges to
$$\frac{\int_X\omega_1^n}{\int_X\omega_0^m\wedge\omega_X^{n-m}}\int_X\eta(\omega_0+\mn\de\db\hat{\varphi})^m\wedge \omega_X^{n-m},$$
and using \eqref{to5}, \eqref{to6} we get
$$\int_X\eta F \omega_{SF}^{n-m}\wedge\omega_0^m=
\frac{\int_X\omega_1^n}{\int_X\omega_0^m\wedge\omega_X^{n-m}}\int_X\eta(\omega_0+\mn\de\db\hat{\varphi})^m\wedge \omega_X^{n-m}.$$
We then integrate first along the fibers and get
\begin{equation*}\begin{split}
\int_{Y}&\eta F \omega_Y^m \left(\int_{X_y}\omega_{SF,y}^{n-m}\right)\\
&=\frac{\int_X\omega_1^n}{\int_X\omega_0^m\wedge\omega_X^{n-m}}\int_Y\eta(\omega_Y+\mn\de\db\hat{\varphi})^m
\left(\int_{X_y} \omega_y^{n-m}\right),
\end{split}\end{equation*}
and since $\omega_y$ is cohomologous to $\omega_{SF,y}$, we get
$$\int_{Y}\eta F \omega_Y^m=
\frac{\int_X\omega_1^n}{\int_X\omega_0^m\wedge\omega_X^{n-m}}\int_Y\eta(\omega_Y+\mn\de\db\hat{\varphi})^m,$$
which is just the weak form of \eqref{wp}. This shows that any weak limit $\hat{\omega}$ of $\ti{\omega}_t$ as $t\to 0$ satisfies \eqref{wp} weakly, and we have already remarked that we can write $\hat{\omega}=\omega_Y+\mn\de\db\hat{\varphi}$ with $\hat{\varphi}$ in $L^\infty$.
By Ko\l odziej's uniqueness of $L^\infty$ weak solutions of \eqref{wp} (see \cite[Theorem 3.2]{SoT2} and \cite{egz1, zhangthesis}), we must have $\hat{\varphi}=\psi$, and so the whole sequence $\ti{\omega}_t$ converges weakly to $\omega$ as $t\to 0$. Then the bound \eqref{c2a} implies that $\varphi_t$ actually converges to $\psi$ in the $C^{1,\beta}_{loc}$ topology on $X\backslash S$.
\end{proof}

\section{Examples and remarks}\label{examp}
In this section we give some examples where Theorem \ref{main2} applies.\\

The easiest example is a complex torus $X$ of dimension $n$ fibering over another torus $Y$ of lower dimension $m$. The fibers are also tori and they are all biholomorphic.
In this case Ricci-flat metrics are just flat, and they can be identified with constant positive definite Hermitian $n\times n$ matrices. 
If we degenerate the K\"ahler class on $X$ to the pullback of a K\"ahler class from $Y$, the matrices converge to a nonnegative definite matrix whose kernel generates the tangent space to the fibers.
So the fibers are shrunk to points and the flat metrics on $X$ converge to the flat metric on $Y$ in the given class. This is of course compatible with Theorem \ref{main2}, because in this case the Weil-Petersson metric is identically zero,
and the set $S$ of singular fibers is empty.

To see a more interesting example, let $X$ be an elliptically fibered $K3$ surface, so $X$ comes equipped with a morphism $f:X\to\mathbb{P}^1$
with generic fibers elliptic curves. Then the pullback of an ample line bundle on $\mathbb{P}^1$ gives a nef line bundle $L$ on
$X$ with Iitaka dimension $1$. In the case when all the singular fibers of $f$ are of Kodaira type $I_1$, Gross-Wilson have
shown in \cite{gw} that sequences of Ricci-flat metrics on $X$ whose class approaches $c_1(L)$ converge in $C^\infty$
on compact sets of the complement of the singular fibers to the pullback of a K\"ahler metric on $\mathbb{P}^1$ (minus the $24$ points which correspond to the singular fibers). Their argument relies on explicit model metrics that are almost Ricci-flat, and it is not well-suited to generalization to higher dimensions.
More recently Song-Tian \cite{st} gave a more
direct proof of the result of Gross-Wilson and they noticed that the limit metric has Ricci curvature equal to the Weil-Petersson metric. Our Theorem \ref{main2} applies in this example, as well as in higher dimensions.

One can easily construct examples of higher-dimensional Calabi-Yau manifolds that are algebraic fiber spaces, to which Theorem \ref{main2} applies. For example the case of Calabi-Yau threefolds is studied extensively in \cite{og}, where many examples are given.



\begin{thebibliography}{2}
\bibitem[A]{allard} Allard, W.K. {\em On the first variation of a varifold},
Ann. of Math. (2) {\bf 95} (1972), 417--491, MR0307015, Zbl 0252.49028.
\bibitem[DP]{dpali} Demailly, J.-P., Pali, N. \emph{Degenerate complex Monge-Amp\`ere equations over compact K\"ahler manifolds}, preprint, arXiv:0710.5109.
\bibitem[EGZ1]{egz1} Eyssidieux, P., Guedj, V., Zeriahi, A. \emph{Singular K\"ahler-Einstein metrics}, J. Amer. Math. Soc. {\bf 22}  (2009), 607-639, MR2505296.
\bibitem[EGZ2]{egz2} Eyssidieux, P., Guedj, V., Zeriahi, A. \emph{A priori $L^\infty$-estimates for degenerate complex Monge-Amp\`ere equations}, Int. Math. Res. Not. {\bf 2008}, Art. ID rnn 070, 8 pp, MR2439574, Zbl 1162.32020.
\bibitem[Fi]{fine} Fine, J. \emph{Constant scalar curvature K\"ahler metrics on fibred complex surfaces}, J. Differential Geom. {\bf 68} (2004), no. 3, 397--432, 
MR2144537, Zbl 1085.53064.
\bibitem[FS]{fs} Fujiki, A., Schumacher, G. {\em The moduli space of extremal compact K\"ahler manifolds and generalized Weil-Petersson metrics},  Publ. Res. Inst. Math. Sci.  26  (1990),  no. 1, 101--183, MR1053910, Zbl 0714.32007.
\bibitem[GW]{gw} Gross, M., Wilson, P.M.H. \emph{Large complex structure limits of $K3$ surfaces}, J. Differential Geom. {\bf 55} (2000), no. 3, 475--546, 
MR1863732, Zbl 1027.32021.
\bibitem[H]{hebey} Hebey, E. {\em Sobolev spaces on Riemannian manifolds},
Lecture Notes in Mathematics, 1635, Springer-Verlag, Berlin, 1996, 
MR1481970, Zbl 0866.58068.
\bibitem[K]{kol} Ko\l odziej, S. \emph{The complex Monge-Amp\`ere equation}, Acta Math. {\bf 180} (1998), no. 1, 69--117, MR1618325, Zbl 0913.35043.
\bibitem[KT]{koltian} Ko\l odziej, S., Tian, G. \emph{A uniform $L\sp \infty$ estimate for complex Monge-Amp\`ere equations}, Math. Ann. {\bf 342} (2008),  no. 4, 773--787, MR2443763, Zbl 1159.32022.
\bibitem[La]{laza} Lazarsfeld, R. \emph{Positivity in algebraic geometry I}, Springer-Verlag, Berlin, 2004, MR2095471, Zbl 1093.14501.
\bibitem[LY]{LY} Li, P., Yau, S.-T. {\em Estimates of eigenvalues of a compact Riemannian manifold}, in {\em Geometry of the Laplace operator (Honolulu, Hawaii, 1979)}, Proc. Sympos. Pure Math. {\bf 36}, Amer. Math. Soc., Providence, 1980, 205--239, MR0573435, Zbl 0441.58014.
\bibitem[McM]{ctm} McMullen, C.T. \emph{Dynamics on $K3$ surfaces: Salem numbers and Siegel disks}, J. Reine Angew. Math. {\bf 545} (2002), 201--233, MR1896103, 
Zbl 1054.37026.
\bibitem[MS]{MS} Michael, J.H., Simon, L.M. \emph{Sobolev and mean-value inequalities on generalized submanifolds of $R\sp{n}$}, Comm. Pure Appl. Math. {\bf 26} (1973), 361--379, MR0344978, Zbl 0256.53006.
\bibitem[O]{og} Oguiso, K. {\em On algebraic fiber space structures on a Calabi-Yau $3$-fold}, Internat. J. Math. {\bf 4} (1993), no. 3, 439--465, MR1228584, Zbl 0793.14030.
\bibitem[PSS]{PSS} Phong, D.H., \v Se\v sum, N., Sturm, J. {\em Multiplier ideal sheaves and the K\"ahler-Ricci flow}, Comm. Anal. Geom. {\bf 15} (2007), no. 3, 613--632, MR2379807, Zbl 1143.53064.
\bibitem[ST1]{st} Song, J., Tian, G. \emph{The K\"ahler-Ricci flow on surfaces of positive Kodaira dimension}, Invent. Math. {\bf 170} (2007), no. 3, 609--653, MR2357504, Zbl 1134.53040.
\bibitem[ST2]{SoT2} Song, J., Tian, G. {\em Canonical measures and K\"ahler-Ricci flow}, preprint, arXiv:0802.2570.
\bibitem[ST3]{SoT3} Song, J., Tian, G. {\em Convergence of the K\"ahler-Ricci flow on minimal models}, in preparation.
\bibitem[St]{stoppa} Stoppa, J. {\em Twisted cscK metrics and K\"ahler slope stability}, preprint, arXiv:0804.0414.
\bibitem[Tp]{topping} Topping, P. {\em Relating diameter and mean curvature for submanifolds of Euclidean space},  Comment. Math. Helv.  {\bf 83}  (2008),  no. 3, 539--546, MR2410779, Zbl 1154.53007.
\bibitem[To1]{schwarz} Tosatti, V. \emph{A general Schwarz Lemma for almost-Hermitian manifolds}, Comm. Anal. Geom. {\bf 15} (2007), no.5, 1063-1086, MR2403195, Zbl 1145.53019.
\bibitem[To2]{deg} Tosatti, V. \emph{Limits of Calabi-Yau metrics when the K\"ahler class degenerates}, J. Eur. Math. Soc. (JEMS) {\bf 11} (2009), no.4, 755-776, Zbl pre05608944.
\bibitem[To3]{tesi} Tosatti, V. {\em Geometry of complex Monge-Amp\`ere equations}, PhD thesis, Harvard University, 2009.
\bibitem[TWY]{taming} Tosatti, V., Weinkove, B., Yau, S.-T. \emph{Taming symplectic forms and the Calabi-Yau equation}, Proc. London Math. Soc. {\bf 97} (2008), no.2, 401-424, MR2439667, Zbl 1153.53054.
\bibitem[W]{wilson} Wilson, P.M.H. \emph{Metric limits of Calabi-Yau manifolds}, in \emph{The Fano Conference}, 793--804, Univ. Torino, Turin, 2004, MR2112603, Zbl 1061.32019.
\bibitem[Y1]{Ya} Yau, S.-T. \emph{On the Ricci curvature of a compact K\"ahler manifold and the complex Monge-Amp\`ere equation, I}, Comm. Pure Appl. Math. {\bf 31} (1978), 339--411, MR0480350, Zbl 0369.53059.
\bibitem[Y2]{yauschw} Yau, S.-T. \emph{A general Schwarz lemma for K\"ahler manifolds}, Amer. J. Math. {\bf 100} (1978), no. 1, 197--203, MR0486659, Zbl 0424.53040.
\bibitem[Y3]{yau2} Yau, S.-T. \emph{Problem section} in \emph{Seminar on Differential Geometry}, pp. 669--706, 
Ann. of Math. Stud. {\bf 102}, Princeton Univ. Press, 1982 (problem 49), MR0645762, Zbl 0479.53001.
\bibitem[Z]{zhangthesis} Zhang, Z. \emph{Degenerate Monge-Amp\`ere equations over projective manifolds}, Ph.D. Thesis, MIT, 2006.
\end{thebibliography}
\end{document}